\documentclass{amsart}
\usepackage{amscd}
\numberwithin{equation}{section}
\let\cal\mathcal

\newcommand{\Nor}{\mathop{\mathrm{Nor}}}

\newcommand{\Span}{\mathop{\mathrm{Span}}}

\newtheorem{lemma}{Lemma}[section]
\newtheorem{proposition}[lemma]{Proposition}
\newtheorem{theorem}[lemma]{Theorem}

\newtheorem{maintheorem}[lemma]{Main Theorem}

\newtheorem{corollary}[lemma]{Corollary}

\newtheorem{facts}[lemma]{Facts}

\theoremstyle{definition}

\newtheorem{example}[lemma]{Example}
\newtheorem{definition}[lemma]{Definition}

\newtheorem{notation}[lemma]{Notation}

{

\newtheorem*{proofmaintheorem}{Proof of the Main Theorem}

\newtheorem*{prooftheorem}{Proof of Theorem 3.1}

}

\theoremstyle{remark}

\newtheorem{remark}[lemma]{Remark}

\newcommand{\Acal}{\mbox{$\cal A$}}

\newcommand{\Dcal}{\mbox{$\cal D$}}

\newcommand{\Lcal}{\mbox{$\cal L$}}

\newcommand{\Ncal}{\mbox{$\cal N$}}

\newcommand{\Ocal}{\mbox{$\cal O$}}

\newcommand{\Wcal}{\mbox{$\cal W$}}
\newcommand{\Rcal}{\mbox{$\cal R$}}

\newcommand{\Scal}{\mbox{$\cal S$}}

\newdimen\uboxsep \uboxsep=1ex
\def\uboxn#1{\vtop to 0pt{\hrule height 0pt depth 0pt\vskip\uboxsep
\hbox to 0pt{\hss #1\hss}\vss}}

\def\uboxs#1{\vbox to 0pt{\vss\hbox to 0pt{\hss #1\hss}
\vskip\uboxsep\hrule height 0pt depth 0pt}}

\title[Garside monoids with quadratic relations]{Garside structure on monoids with
quadratic square-free relations} \keywords{Yang-Baxter, Garside
monoids, Garside groups, Binomial skew-polynomial semigroups,
Frobenius Algebras, Quadratic algebras,  Quantum Groups}
\subjclass{Primary 81R50, 16S37, 16W35, 16W50, 16S36, 81R60}
\thanks{The author was partially supported by the
 The Abdus Salam International Centre
   for Theoretical Physics (ICTP), Trieste}
\author{Tatiana Gateva-Ivanova}
\address{Institute of Mathematics and Informatics\\
Bulgarian Academy of Sciences\\
Sofia 1113, Bulgaria\\
and
American University in Bulgaria\\
2700 Blagoevgrad, Bulgaria }

\email{tatianagateva@yahoo.com, tatyana@aubg.bg}

\begin{document}
\date{\today}
\begin{abstract}
We show the intimate connection between various mathematical notions
 that are currently under active investigation: a class of \emph{Garside monoids},
with a ``nice" Garside element, certain monoids $S$ with quadratic
relations, whose monoidal algebra $A= \textbf{k}S$ has a
\emph{Frobenius Koszul dual} $A^{!}$ with \emph{regular socle},
the monoids of \emph{skew-polynomial type} (or equivalently,
\emph{binomial skew-polynomial rings}) which were introduced and
studied by the author and in 1995 provided a new class of
Noetherian \emph{Artin-Schelter regular}
 domains,
and \emph{the  square-free set-theoretic  solutions of the Yang-Baxter equation}.
There is a beautiful symmetry in these objects due to their nice combinatorial and algebraic properties.
\end{abstract}
\maketitle
\section{Introduction}

Let $X$ be a
nonempty set and let $r: X\times X \longrightarrow X\times X$ be a
bijective map. In this case  we shall use notation $(X,r)$ and
refer to it as a \emph{quadratic set} or set with quadratic map $r$. We present the
image of $(x,y)$ under $r$ as
\begin{equation}
\label{r} r(x,y)=({}^xy,x^{y}).
\end{equation}
The formula (\ref{r}) defines a ``left action" $\Lcal: X\times X
\longrightarrow X,$ and a ``right action" $\Rcal: X\times X
\longrightarrow X,$ on $X$ as:
\begin{equation}
\label{LcalRcal} \Lcal_x(y)={}^xy, \quad \Rcal_y(x)= x^{y},
\end{equation}
for all $ x, y \in X.$ The map  $r$ is \emph{left nondegenerate (respectively, right nondegenerate)}, if the maps $\Lcal_x$ (respectively,  $\Rcal_x$) are
bijective for each $x\in X$. We say that $r$ is {\em nondegenerate} if it is both left and right nondegenerate,
$r$ is \emph{involutive} if $r^2 = id_{X \times X}$.  In this paper we shall  always  assume that $r$ is nondegenerate. Also, as a notational tool, we  shall often
identify the sets $X^{\times k}$ of ordered $k$-tuples, $k \geq 2,$  and $X^k,$ the set of
all monomials of length $k$ in the free monoid $\langle
X\rangle.$

As in \cite{T04, TSh07, TSh08, TSh0806, TP}, to each quadratic set
$(X,r)$ we associate canonically algebraic objects (see
Definition~\ref{associatedobjects}) generated by $X$ and with
quadratic defining relations $\Re =\Re(r)$ naturally determined as
\begin{equation}
\label{defrelations}
\begin{array}{ll}
xy=y^{\prime} x^{\prime}\in \Re(r),&\text{whenever} \\
 r(x,y) = (y^{\prime}, x^{\prime}) & \text{and} \quad (x,y) \neq (y^{\prime}, x^{\prime})\;\;\text{hold in}\;\;X \times X.
\end{array}
\end{equation}
Note that in the case when $X$ is finite, the set $\Re(r)$  of
defining
 relations is also finite, therefore the associated algebraic
 objects are finitely presented. We can tell the precise number of defining relations,  whenever $(X,r)$  is nondegenerate and involutive, with $\mid X \mid =n$. In this case the set of defining relations $\Re(r)$,  contains precisely $\binom{n}{2}$ quadratic
 relations (see  Proposition \ref{nondeg_and_invol_goodprop}  ).
 In many cases the associated algebraic
 objects  will be \emph{standard finitely presented} with respect to the
 degree-lexicographic ordering induced by an appropriate enumeration of $X$.  It is known in particular that the algebra generated by the monoid
 $S(X,r)$ defined in this way  has remarkable homological properties when $r$ is involutive, nondegenerate
 and obeys the braid or Yang-Baxter equation
 in $X\times X\times X$. Set-theoretic solutions were introduced in  \cite{D,W} and have been under intensive study during the last decade.
 There are many  works on set-theoretic solutions and related
structures, of which a relevant selection for the interested reader
is \cite{W, TM, RTF,  ESS, Lu,  T04, T04s, Rump, Takeuchi, Veselov,
TSh07, TSh08, TSh0806, TP, Ch}. In this section we recall basic
notions and results which will be used in the paper.
 We shall use the terminology,  notation and some results from
 our previous works,
\cite{TLovetch, T04s,  T04, TSh08}. This paper is a natural continuation of \cite{T04s} which provides results of significant importance for the Main Theorem.

  Here we study the  close relations between various notions that appear in the literature:
 square-free solutions of the Yang-Baxter equation;  a class of Garside structures (we call them \emph{regular Garside monoids});  skew polynomial monoids;
 semigroups of I type; and a class of monoids $S$ whose monoidal algebra $A=kS$ has
 Frobenius Koszul dual $A^{!}$ with special properties. All these are also closely related
 to Artin-Schelter regularity of certain quadratic algebras.

For a non-empty set $X$, as usual, we denote by $\langle X \rangle$,
and ${}_{gr}\langle X \rangle$, respectively, the free  monoid and
the free group generated by $X$ and by $\textbf{k}\langle X
\rangle$- the free associative $\textbf{k}$-algebra generated by
$X$, where $\textbf{k}$ is an arbitrary field. For a set $F
\subseteq \textbf{k}\langle X \rangle$, $(F)$ denotes the two sided
ideal of $\textbf{k}\langle X \rangle$ generated by $F$.

\begin{definition} \cite{T04, TSh08}
\label{associatedobjects} Assume that $r:X^2 \longrightarrow X^2$
is a bijective map.

(i) The monoid
\[
S =S(X, r) = \langle X ; \; \Re(r) \rangle,
\]
 with a
set of generators $X$ and a set of defining relations $ \Re(r),$
is called \emph{the monoid associated with $(X, r)$}.

(ii) The \emph{group $G=G(X, r)$ associated with} $(X, r)$ is
defined as
\[
G=G(X, r)={}_{gr} \langle X ; \; \Re (r) \rangle.
\]

(iii) For arbitrary fixed field $\textbf{k}$, \emph{the
$\textbf{k}$-algebra associated with} $(X ,r)$ is defined as
\begin{equation}
\label{Adef} \Acal = \Acal(\textbf{k},X,r) = \textbf{k}\langle X ; \;\Re(r) \rangle \simeq \textbf{k}\langle X  \rangle /(\Re_0),
\end{equation}
where $\Re_0 = \Re_0(r)$ is the set of \emph{binomials} in $k\langle X  \rangle $:
\begin{equation}
\label{Re0} \Re_0 = \{xy-y^{\prime}x^{\prime}\mid xy=y^{\prime}x^{\prime}\in \Re (r) \}.
\end{equation}
\end{definition}
Clearly $\Acal$ is a quadratic algebra, generated by $X$ and
 with defining relations  $\Re_0(r).$
Furthermore, $\Acal$  is isomorphic to the monoid algebra
$\textbf{k}S(X, r).$

\begin{remark}
\label{quadraticmonoid/quadraticsetremark} Conversely, as in
\cite{TM, T04s}, each finitely presented monoid $S =\langle X ; \;
\Re \rangle$, where $\Re$ is \emph{a set of quadratic binomial
relations} $xy=y^{\prime} x^{\prime}$ such that each monomial $xy\in
X^2$ occurs at most once in $\Re $ defines canonically a quadratic
set $(X, r)$. Let $r$ be the involutive bijective map determined via
$\Re $ as follows
\[
r(x,y) = \left\{ \begin{array}{l} (y^{\prime}, x^{\prime}),\; \text{if}\;
xy=y^{\prime}x^{\prime}\in \Re  \\
                                       (x,y)\quad \text{else}.
\end{array}\right\}
\]
Clearly, $(X, r)$ is a quadratic set and $S \simeq S(X, r).$
\end{remark}

\begin{remark}
\label{orbitsinG}
When we study the monoid $S=S(X,r)$, the group $G=G(X,r)$,  or
the algebra $\Acal = \Acal(\textbf{k}, X,r) \simeq \textbf{k }[ S]$ associated with $(X,r)$,
it is  convenient to use the action of the infinite groups,
$\Dcal_k(r)$, generated by maps associated with the quadratic
relations, as follows. We consider the  bijective
maps
\[
r^{ii+1}: X^k \longrightarrow X^k, \; 1 \leq i \leq k-1, \quad\text{where}\quad r^{ii+1}=Id_{X^{i-1}}\times r\times Id_{X^{k-i-1}}.
\]
Note that these maps are elements of the symmetric group $Sym
(X^k)$.
  Then the group $\Dcal _k(r)$ generated by  $r^{ii+1}, \; 1 \leq i \leq k-1,$
acts on $X^k.$ If  $r$ is involutive, the bijective maps $r^{ii+1}$
are involutive, as well, so in this case $\Dcal_k(r)$ is
the infinite group
\begin{equation}
\label{Dk}
 \Dcal_k(r)=\: _{\rm{gr}} \langle r^{ii+1}\mid\quad (r^{ii+1})^2= e,\quad  1 \leq i \leq k-1 \rangle.
 \end{equation}

Note first that $\Dcal_k(r)$ is isomorphic to $(C_2)^{*(k-1)},$
the free product of $k-1$ cyclic groups of order 2. In particular,
for $k=3$,  $\Dcal_3(r)$ is simply the \emph{the infinite dihedral
group}.

Secondly, note that if $\omega$ is a monomial of length $k$ in $\langle X\rangle$,
the set of all monomials $\omega^{\prime}\in \langle X\rangle$ such that $\omega^{\prime} =\omega$ as elements of $S,$ coinside with
the orbit    $\Ocal_{\Dcal_k(r)}(\omega)$ in $X^k.$ Analogous statements are true for  $G$ and $\Acal$.
\end{remark}
\begin{definition}
\label{quantumbinomialset&ybedef} Let $(X,r)$ be a quadratic set.
\begin{enumerate}
 \item
$(X,r)$ is said to be \emph{square-free} if $r(x,x)=(x,x)$ for all $x\in X.$
\item  $(X,r)$ is called a \emph{quantum binomial
set} if it is nondegenerate, involutive and square-free.
\item \label{YBE} $(X,r)$ is a \emph{set-theoretic solution of the
Yang-Baxter equation}
 (YBE)   if  the braid relation
\[r^{12}r^{23}r^{12} = r^{23}r^{12}r^{23}\]
holds in $X\times X\times X.$  In this case $(X,r)$ is also called a \emph{braided set}. If in addition $r$ is involutive $(X,r)$ is called
\emph{a symmetric set}.
\end{enumerate}
\end{definition}

 In various cases in a nodegenerate $(X,r)$ the left and the
right actions are inverses, i.e. $\Rcal_x= \Lcal_x^{-1}$ and
$\Lcal_x = \Rcal_x^{-1}$ for all $x\in X$. For example
this is true for every square-free symmetric set $(X,r)$, see \cite{TLovetch, T04}. In   \cite{TSh08},
we singled out a class of nondegenerate sets
$(X,r)$ by a
condition \textbf{lri} defined below.

\begin{definition} \cite{TSh08}
\label{lri}
Let  $(X,r)$ be a quadratic set (of arbitrary cardinality). We define the condition
\begin{equation}
\label{lrieq}
\textbf{lri:} \quad\quad ({}^xy)^x= y={}^x{(y^x)} \;\text{for all} \; x,y \in X.
\end{equation}
In other words \textbf{lri} holds if and only if  $\Rcal_x=\Lcal_x^{-1}$ and $\Lcal_x = \Rcal_x^{-1}$.
\end{definition}

\begin{definition} \cite{KasselT}
\label{garsidedef}
\begin{enumerate}
\item
\label{pregarsidemonoiddef}
\emph{A pre-Garside monoid } is a pair consisting of a monoid $M$ and an
element $\Delta\in M$
such that the set $\Sigma= \Sigma_{\Delta}$ of all left divisors of $\Delta$
satisfies the following conditions:

(a) $\Sigma$ is finite, generates $M$, and coinsides with the set of right divisors of $\Delta.$

(b) If $a,b \in \Sigma$  are such that $\Delta a = \Delta b,$ or $a\Delta= b\Delta,$ then $a=b$

The element $\Delta \in M$ is called \emph{the Garside element of $M$}.
\item
\label{garsidemonoiddef} Let $(M, \Delta)$ be a pre-Garside monoid
and let $\Sigma$ be the set of left (and right) divisors of
$\Delta.$ The pair $(M, \Delta)$ is \emph{a Garside monoid } if $M$
is atomic and for any two atoms $s,t$ of $M$, the set
\[
\{a \in \Sigma \mid s \preceq a, \;\text{and}\; t \preceq a\}
\]
has a (unique) minimal element $\Delta_{s,t}$ (with respect to $\preceq$).

Here  $\preceq$ is the partial order on $M$ defined as $s \preceq a,$ if $s$ is a left divisor of $a$ that is $a = s\alpha$ for some $\alpha \in M$.
It is known that each Garside monoid $M$ embeds in its group of fraction $G_M.$ $G_M$ is called \emph{a Garside group}.
\item
\label{comprehensivegarsidemonoiddef}
A set $\Sigma \subset M$ is \emph{comprehensive  } if $1 \in \Sigma$  and $M$ has a presentation by generators and relations such that all generators and relators belong to $\Sigma.$

A Garside monoid $(M,\Delta)$ is\emph{ comprehensive } if the set $\Sigma$  of the divisors of $\Delta$
is comprehensive.
\end{enumerate}
\end{definition}
\begin{remark}
\label{atomicremark} By definition $S= \langle X; \Re \rangle$ has
quadratic defining relations, thus, it has \emph{a length balanced
presentation } so by \cite[Lemma 6.4]{KasselT},  $S(X,r)$ is atomic,
and the set of its atoms coincides with $X$.
\end{remark}

We introduce now the notions of \emph{regular quantum monoids} and
\emph{regular Garside monoids}. In each case the name
\emph{regular } comes in a natural way-we show that the two
notions are equivalent, and the monoids have a rich list of 'good'
properties. Furthermore, in each case the semigroup algebra
$\textbf{k}S$ is \emph{Artin-Schelter regular ring}, see Corollary
\ref{ASregularitycor}.   Artin-Schelter regularity, see \cite{AS},
is  an important notion in noncommutative geometry.

\begin{definition}
\label{regquantummondefinitions}
Let $(X,r)$ be a finite quantum binomial set with \textbf{lri}, and let $\mid X \mid =n.$
Let
$S= S(X,r),
\Acal = \Acal(\textbf{k},X,r)$ be respectively the associated monoid,  and the associated quadratic algebra.
Let $\Acal^{!}$
be the Koszul dual of $A$, see Definition \ref{koszuldualdef}.
\begin{enumerate}
\item
\label{regularmonomialdef1}
A monomial $\omega \in S$ of length $k$ is called \emph{a square-free monomial}, if its orbit
$\Ocal_{\Dcal_k(r)}(\omega)$ in $X^k$
does not contain words of the shape $\omega^{\prime} = a xx b,$ where $x \in X, a, b \in \langle X \rangle$.
In other words there is no equality $\omega = axxb$ as elements of $S$.
\item
Let $\omega_0 = x_1x_2 \cdots x_n \in X^n$ be a square free element of $S$, such that all $x_i$ are pairwise distinct, so we
enumerate $X = \{ x_1, \cdots, x_n \}$ and fix the degree-lexicographic ordering $ < $ on $\langle X \rangle$, where $x_1 < x_2 < \cdots < x_n$. We say that $\omega_0$ is a regular element of $S$  if it is the minimal element, with respect to $ < $, in the orbit $\Ocal_{\Dcal_k(r)}(\omega_0).$  We shall also say that each $\omega \in \Ocal_{\Dcal_k(r)}(\omega_0)$  has \emph{a regular presentation} $\omega=\omega_0$.

 \item
\label{regquantummondef2}
$S$ is called \emph{a regular quantum monoid } if

 (i) The Koszul dual $\Acal^{!}$ is Frobenius  of dimension $n$, and

 (ii) The \emph{principal monomial} $W$ of $S$ has \emph{a regular presentation} $W = x_1x_2 \cdots x_n.$

 ( By definition $W$ spans the socle of $\Acal^{!}$ and is the longest square-free element in $S$, see section \ref{section_regularity}).
 \item
 \label{regquantummondef3}
 A Garside monoid $S$ is \emph{a regular Garside  monoid} if

  (i) $S$ has a Garside element $\Delta \in S$ with a regular presentation $\Delta = x_1x_2...x_n$
 and

 (ii) every square-fee monomial $a \in S$ of length $\leq n$ is a left (and a right) divisor of $\Delta.$

 In this case $\Delta$ is called  \emph{a regular Garside element of $S$}. Clearly, $\Delta$ is unique.
 \end{enumerate}
\end{definition}
\begin{remark}
Note that, in general, a regular element $\omega \in S$ can have more than one regular presentations, see Example \ref{moreregpresentationsex}.
We don't know whether it is  possible that  the  monoid $S(X,r)$ of an arbitrary quantum binomial set $(X,r)$ with \textbf{lri}, may have two distinct regular elements $\omega, \omega^{\prime}$ (in the sence that $\omega \neq \omega^{\prime}$ as elements of $S$).
 \end{remark}
\begin{example}***
\label{moreregpresentationsex}
\end{example}
\begin{definition} \cite{T96}
\label{skewpolynomialsemigroup} We say that the monoid $S$ is
 \emph{
a monoid of skew-polynomial type}, (or shortly, \emph{a
skew-polynomial semigroup}) if it has a standard finite
presentation as
\[
S  = \langle X; \Re \rangle, \] where the set of generators $X$ is
ordered: $x_1 < x_2 < \cdots < x_n,$ and $\Re$ is a set of
$\binom{n}{2}$ quadratic
 relations,
\begin{equation}
\label{skewrelations}
\Re  =\{ x_jx_i=x_{i^{\prime}}x_{j^{\prime}})\mid 1 \leq i< j
\leq n, 1 \leq i^{\prime} < j^{\prime} \leq n \},
\end{equation}
satisfying
\begin{enumerate}
\item[(i)] \label{skewpolynomialsemigroupi} each monomial $xy\in
X^2$, with $x\neq y$, occurs in exactly one relation in $\Re $ (a
monomial of the type $xx$ does not occur in any relation in
$\Re$); \item[(ii)] \label{skewpolynomialsemigroupii} if
$(x_jx_i=x_{i^{\prime}}x_{j^{\prime}})\in \Re $, with $1 \leq i< j
\leq n,$ then  $i^{\prime} < j^{\prime}$, and $j > i^{\prime}$.
(this also imply $i < j^{\prime}$, see \cite{T96}) \item[(iii)]
\label{skewpolynomialsemigroupiii} the overlaps $x_kx_jx_i$ with
$k>j>i, \; 1\leq i, j, k \leq n$ do not give rise to new relations
in $S$, or equivalently,  see \cite{B}, the set of polynomials
$\Re_0$ is a Gr\"{o}bner basis of the ideal $(\Re_0)$ with respect
to the degree-lexicographic ordering of the free semigroup
$\langle X \rangle$.
\end{enumerate}
\end{definition}
\begin{remark}
Clearly, the set of relations (\ref{skewrelations}) defines canonically an involutive   bijective map
\[\begin{array}{c}
r: X\times X \longrightarrow X\times X, \\
  (x_j,x_i) \leftrightarrow (x_{i^{\prime}},x_{j^{\prime}}), \;
  1 \leq i< j \leq n, \;j> i^{\prime}\; 1 \leq i^{\prime} < j^{\prime}\leq n\\
   (x,x) \leftrightarrow (x,x), \forall x \in X.
\end{array}\]
and $(X,r)$ is a quantum binomial quadratic set.
\end{remark}
\begin{remark}
The theory of (noncommutative) Gr\"{o}bner bases implies that each
monomial  $u\in S_0$ has a uniquely determined \emph{normal form},
$\Nor(u)$,  that is the minimal (w.r.t. $<$) element in the orbit
$\Dcal_m(u),$ where $\mid u \mid =m.$ It is known that the normal
form can be found effectively by applying finite steps of reductions
determined via the relations. When the relations are of
skew-polynomial type, as in \ref{skewrelations}, the normal form
$\Nor(u)$ is \emph{an ordered monomial}, ($ x_{i_1} \cdots x_{i_m},$
with $i_1 \leq i_2 \leq \cdots \leq i_m$), and (since $\Re$ is a
Gr\"{o}bner basis) any ordered monomial in $\langle X \rangle$ is in
normal form (mod $(\Re)$). In other words,
condition  (iii) of Definition \ref{skewpolynomialsemigroup} may be rephrased
by saying that  (as sets) we can identify $S_0$ with \emph{the set of ordered monomials}
\[
\label{ncal} \Ncal _0 = \{x_{1}^{\alpha_{1}}\cdots
x_{n}^{\alpha_{n}}\mid \alpha_{n} \geq 0 \text{ for } 1 \leq i
\leq n\}.
\]
\end{remark}

We recall that the skew polynomial monoids appeared first in the
context of a class of quadratic algebras with a
Poincar\'{e}-Birkhoff-Witt type $\textbf{k}$-basis, namely
\emph{the binomial skew polynomial rings}. These rings were
introduced and studied in \cite{T94, T96, T96Preprint, TM, T04s,
Laffaille}. Laffaille calls them \emph{quantum binomial algebras}
and uses computer programme to show in \cite{Laffaille}, that for
$\mid X \mid \leq 6$, the associated automorphism $R$ is a
solution of the Yang-Baxter equation.   We prefer to keep the name
``\emph{binomial skew polynomial rings}" since we have been using
this name for already 15 years. We have verified in two different
ways that the binomial skew polynomial rings provide a new  (at
that time) class of Artin-Schelter regular rings of global
dimension $n$, where $n$ is the number of generators $X.$ In fact,
the  challenge was to prove the Gorenstein condition. The very
first proof (1994) which the author reported on the Artin's
seminar MIT, involved spectral sequences. A second proof (which is
not yet well known) is combinatorial, it involves the Frobenius
property of the Koszul dual algebra $A^{!}$, see
\cite{T96Preprint}. This result is also published in
 \cite{T04s}, Section 3. The third proof which is now well-known, was given in \cite{TM},
 and uses the good homological and algebraic properties of semigroups of I type and their
 semigroup algebras. Furthermore, a close relation
 was found in  \cite{TM} between the three notions of semigroups of skew-polynomial type,
semigroups of I-type and a class of set-theoretic solutions of YBE,
the square-free solutions. In 1996 the author conjectured that the
three notions are equivalent, see \cite{T96Conjecture} (and also
\cite{TLovetch}, \cite{TConstanta}. The verification took several
years, but gradually  we better understood the rich structure and
beautiful symmetry in the associated algebraic objects. These are
now used to verify their  close relation with a class of Garside
structures. The study of Garside monoids and Garside groups has been
recently intensified, Garside structures are used  to prove
recognizability of various properties in groups and have significant
importance for the study of braid groups. Recently Chouraqui showed
that for every finite nondegenerate symmetric set $(X,r)$, the
associated monoid $S(X,r)$ is a Garside monoid, hence the group
$G(X,r)$ is a Garside group.

We give priority here to the natural equivalence of the notions of
various algebraic objects each of which is related to a quadratic
algebra, that is  "\emph{a quantum space}" in the sense of Manin,
\cite{Maninpreprint}, Section 3. Their equivalence makes it
possible to combine and explore  the ``good" properties of these
objects.

\begin{maintheorem}
Let $(X,r)$ be a finite quantum binomial quadratic set with \textbf{lri} (see (\ref{lrieq})).
Let $S=S(X,r)$ be the associated monoid.
The following conditions are equivalent
\begin{enumerate}
\item
\label{garside}
$S$ is a regular Garside monoid. \item \label{regularmonoid} $S$
is a regular quantum monoid. \item \label{skewpolymonoid} $S$ is a
monoid of skew-polynomial type, (with respect to some appropriate
enumeration of $X$), or equivalently, $A(k, X,r)$ is a binomial
skew-polynomial ring. \item \label{YBEcond} $(X,r)$ is a
set-theoretic solution of the YBE.
\end{enumerate}
Each of these conditions implies that $(S, \Delta)$ is a
comprehensive Garside monoid with a Garside element $\Delta$ and
its group of fractions is isomorphic to the group $G= G(X,r)$
(with the same generators and relations). In particular
$(G,\Delta)$, is a Garside group.
\end{maintheorem}

\begin{corollary}
\label{ASregularitycor} Let $(X,r)$ be a finite quantum binomial
quadratic set with \textbf{lri}. Each of the conditions (1), (2),
(3), (4)  of the Main Theorem implies that the associated
quadratic algebra $\Acal = \Acal(\textbf{k}, X, r)$ is an
Artin-Shelter regular ring of global dimension $n,$ where $n$ is
the cardinality of $X$. Furthermore, $\Acal$ is Koszul and left
and right Noetherian domain.
\end{corollary}

\
\section{Preliminaries on quadratic sets $(X,r)$ and the associated objects}


In this section $(X,r)$ is a nondegenerate quadratic set of
arbitrary cardinality. In the cases when $X$ is finite this will be
clearly indicated. When this is possible we shall avoid a fixed
enumeration of $X$, to make most of the properties invariant on the
enumeration (see for example the cyclic conditions).

\begin{remark}
\label{someproperties} Suppose  $(X,r)$ is a quadratic set, and let
${}^x\bullet$, and ${\bullet}^x$ be the associated left and right
actions. Then
\begin{enumerate}
\item \label{involutive} $r$ is involutive if and only if
\[
{}^{{}^xy}{(x^y)} = x, \; \text{and}\;
({}^yx)^{y^x} = x, \; \text{for all}\; x,y \in X.\]
 \item \label{SF} $r$ is square-free if and only if
\[
{}^xx=x,\; \text{and} \; x^x = x \; \text{for all}\; x \in X.
\]
\item \label{SFnondeg} If $r$ is nondegenerate and square-free,
then
\[
{}^xy=x\ \Longleftrightarrow\ x^y = y  \Longleftrightarrow\  y=x\ \Longleftrightarrow\ r(x,y)=(x,y).\]
\end{enumerate}
\end{remark}

It is also straightforward to write out the Yang-Baxter equation for
$r$ in terms of the actions. This is in  \cite{ESS} (see also
\cite{TSh08}),  but we recall it here in our notations for
convenience.
\begin{remark}\label{ybe} Let $(X,r)$ be given in the notations above. Then $r$ obeys the YBE
(or $(X,r)$ is a braided set) {\em iff} the following conditions
hold
\[
\begin{array}{lclc}
 {\bf l1:}\quad& {}^x{({}^yz)}={}^{{}^xy}{({}^{x^y}{z})},
 \quad\quad\quad
 & {\bf r1:}\quad&
{(x^y)}^z=(x^{{}^yz})^{y^z},
\end{array}\]
 \[ {\rm\bf lr3:} \quad
{({}^xy)}^{({}^{x^y}{z})} \ = \ {}^{(x^{{}^yz})}{(y^z)},\]
 for all $x,y,z \in X$.
 \end{remark}

Some of the conditions of the following proposition can be
extracted from \cite{TSh08}, and possibly from other sources, see
for example \cite[Lemma 2.12, Corollary 2.13]{TSh08}. However, we
prefer to present here a compact easy proof.
\begin{proposition}
\label{nondeg_and_invol_goodprop}
Suppose $(X,r)$ is a quadratic set of arbitrary cardinality, and   $r$ is nondegenerate and involutive.
Then the following conditions hold.
\begin{enumerate}
\item
$S$ satisfies cancelation law on monomials of length $2$,
that is  for every $s,a,b \in X$
the following implications hold
\begin{equation}
\label{2-cancelation}
\begin{array}{ll}
(a) \quad& sa=sb \;\text{holds in} \; S \Longrightarrow a=b\\
(b) \quad& as=bs \;\text{holds in} \; S \Longrightarrow a=b\\
\end{array}
\end{equation}
\item
\textbf{the left Ore condition on $X^2$:} for every pair $s, t\in X$
there exists a unique pair $a,b\in X$, such that $sa=tb$.
Furthermore, $a=b$ \emph{iff} $s = t$

\item
\textbf{the right Ore condition on $X^2$:}
for every pair $s, t\in X$ there exists a unique pair $a,b\in X$, such that $as=bt$.
Furthermore, $a = b$ \emph{iff } $s = t.$
\item
For each $x\in X$ there exists unique $y \in X$ such that $r(x,y) = (x,y).$
Moreover, if $X$ is finite, of order $n$, then the set of defining relations
$\Re(r)$ contains exactly $\binom{n}{2}$ relations, and $r$ has exactly $n$ "fixed points,"
 that is pairs $(x,y)$ with $r(x,y) = (x,y).$
\end{enumerate}
\end{proposition}
\begin{proof}
We shall prove part (1).
Suppose $sa=sb$ holds in $S$. Then (by the involutiveness of $r$) exactly two cases are possible
\[
\begin{array}{ll}
\textbf{(i)} \quad &(s,a)=(s,b) \; \text{holds in}\; X\times X, \\
\textbf{(ii)} \quad& r(s,a)= (s,b)\;\;\text{and}\;\;r(s,b)= (s,a)\;
\text{hold in}\;X\times X.
 \end{array}
\]
In the case \textbf{(i)} the equality $a = b$ is straightforward.
Assume \textbf{(ii)} holds.
Then
\[
\begin{array}{ll}
 (s,b)= r(s,a)= ({}^sa, s^a)&\quad \Longrightarrow {}^sa= s\\
(s,a)= r(s,b)= ({}^sb, s^b)&\quad \Longrightarrow {}^sb= s
\end{array}\]
This yields
\[
{}^sa= s\quad \text{and}\quad {}^sb= s,
\]
which, by the (left) non degeneracy of $r$ implies $a = b.$

We shall prove next part (2) of the lemma.
Suppose $s,t \in X.$ Then by the nondegeneracy of $r$ there exists unique
$a \in X,$ such that
${}^sa=t.$ Denote $s^a =b.$ Then
clearly, $r(s,a) = ({}^sa, s^a)= (t,b)$ implies
\[
sa = tb \quad\text{holds in}\; S.
\]
That this equality determines $a$ and $b$ uniquely, follows from the \emph{"$2$-cancelation low"} (\ref{2-cancelation}).
The same law implies that  $a=b$ \emph{iff } $s=t$.

Using the right nondegeneracy of $r$ and analogous argument one proves part (3) of the proposiutuion.

We now prove part (4) of the proposition.
Let $x\in X.$
Note first that  for all $ x,y \in X$ the following implication holds.
\begin{equation}
\label{eqf1}
{}^xy = x \Longrightarrow r(x,y) = (x,y)
\end{equation}
Indeed, suppose ${}^xy = x$. Then there are equalities in $X\times X:$
\[
r(x,y) = ({}^xy, x^y) = (x, x^y)
\]
Therefore, there is an equality in $S$
\[
xy = xx^y,
\]
which, by the cancelation law (on generators), see (\ref{2-cancelation})  implies
$x^y = y.$

By the nondegeneracy of $r$, for each $x$ there exists a $y= y_x\in
X$ such that ${}^xy= x,$ and therefore, by (\ref{eqf1}), $r(x,y) =
(x, y)$. We claim that $y=y_x$ is unique with this property. Indeed,
an equality $r(x,z) = (x,z)$ implies ${}^xz= x = {}^xy,$ which by
the (left) nondegeneracy of $r$ implies $z = y.$

Suppose now that $\mid X\mid = n$. Then there are exactly $n$
"fixed" pairs $(x,y_x)$ with the property $r(x,y_x)= (x, y_x).$ The
remaining $n(n-1)$ pairs in $(X\times X$ satisfy $r(a,b) \neq
(a,b)$, which gives exactly one relation
 \[
 ab = {}^ab.a^b \in \Re(r).
 \]
 Clearly, then the number of defining relations is exactly $\binom{n}{2}.$
 This proves part (4). The propsition has been verified.
\end{proof}


***

We recall  the notion of "cyclic conditions" in
terms of the left and right actions

\begin{definition}
\label{cyclicconditionsall} Let $(X,r)$ be a quadratic
set. We define the conditions
\[\begin{array}{lclc}
 {\rm\bf cl1:}\quad&  {}^{y^x}x= {}^yx \quad\text{for all}\; x,y \in X;
 \quad&{\rm\bf cr1:}\quad &x^{{}^xy}= x^y, \quad\text{for all}\; x,y \in
X;\\
 {\rm\bf cl2:}\quad &{}^{{}^xy}x= {}^yx, \quad\text{for all}\; x,y \in X; \quad & {\rm\bf
cr2:}\quad  &x^{y^x}= x^y \quad\text{for all}\; x,y \in X.
\end{array}\]
We say $(X,r)$ is {\em weak cyclic} if {\bf cl1,cr1}
hold and is {\em cyclic} if all four of the above hold.
\end{definition}

One can also define left-cyclic as {\bf cl1,\; cl2} and similarly right-cyclic.

\begin{proposition} \cite{TSh08}
\label{lriprop}
Let  $(X,r)$ be a quadratic set (not necessarily square-free or finite). Then any two of the following conditions imply the remaining third
condition.
\begin{enumerate}
\item \label{cond1} $(X,r)$ is involutive
\item \label{cond2} $(X,r)$ is nondegenerate and cyclic.
\item \label{cond3} {\bf lri} holds.
\end{enumerate}
\end{proposition}

\begin{remark}
\label{lriremark1}
In Proposition~\ref{lriprop} one can replace (2) by the weaker condition\\
\noindent (2')\quad $(X,r)$ nondegenerate and  {\bf cl1} holds.
\end{remark}
\begin{corollary}
\label{lricor}
Every quantum binomial quadratic set $(X,r)$ with \textbf{lri} satisfies all cyclic conditions, see Definition \ref{cyclicconditionsall}.
\end{corollary}

Clearly, \textbf{lri} implies that whatever property is satisfied
by the left action, an analogous property is valid for the right
action and vice versa. In particular, this is valid for the left
and right `cyclic conditions'. The cyclic conditions were
discovered first in 1990, when the author studied
 binomial rings with skew polynomial relation.  It is interesting to know that the proofs of the good algebraic
and homological properties of these algebras and monoids use in
explicit or implicit form the existence of the full cyclic
condition, see Definition~\ref{cyclicconditionsall}. This includes
the properties of being Noetherian, Gorenstein, therefore
Artin-Shelter regular, being of $I$-type and "producing" solutions
of YBE, see \cite{T94,T96,T96Preprint,T04,T04s,TM, GJO, JO}. In
1996 the author was aware that each (finite) square-free solution
$(X,r)$ satisfies the (full) cyclic conditions. This was reported
in various talks as one of the evidences for our conjecture that
every square-free solution of finite order can be "generated" from
a  binomial ring with skew polynomial relation, see for example
\cite{TConstanta}, and \cite{TLovetch}.

Compared with these
works, in \cite{TSh08}  we do not assume that $X$ is finite,  initially the
only restriction on the map $r$ we impose is "$r$ is
nondegenerate".
and study the implication of
the cyclic conditions on the properties of the actions.

\begin{remark}
\label{lriremark2}
The result of these paper extend in a natural way (in direction of Garside structures) the results in \cite{T04s}, which are essential in  various proofs of the new results. Note that in most statements here we assume $(X,r)$ satisfies \textbf{lri}, while in \cite{T04s} we have assumed \emph{the weak cyclic conditions}, see Definition \ref{cyclicconditionsall}. However, these two conditions are equivalent for $(X,r)$ nondegenerate ansd involutive, see Proposition \ref{lriprop} and Remark \ref{lriremark1}. Still \textbf{lri} seems easy and natural to formulate and convenient to use  as it guarantees that the left actions $\Lcal_x, x \in X$ determine uniquely  the right actions ($\Rcal_x = (\Lcal_x)^{-1}$ and therefore $r$. \textbf{lri} is always satisfied in the case when $(X,r)$ is a square-free solution of YBE (of arbitrary cardinality).
\end{remark}

The following can be extracted from \cite[Theorem 2.35]{TSh08},
where more equivalent condtitions are stated.

\begin{theorem}\label{qbs} \cite{TSh08}
Suppose $(X,r)$ is a quantum binomial quadratic set of arbitrary cardinality (i.e. non-degenerate,
involutive and square-free). Then the following conditions are
equivalent:
\begin{enumerate}
\item \label{YBESquarefree} $(X,r)$ is a set-theoretic
solution of the Yang-Baxter equation. \item \label{YBEl1} $(X,r)$
satisfies {\bf l1}.
\item \label{YBEr1} $(X,r)$ satisfies {\bf r1}.
 \item\label{YBElr3}
$(X,r)$ satisfies {\bf lr3}


In this case $(X,r)$ is cyclic and satisfies   \textbf{lri}.
\end{enumerate}
\end{theorem}

Several other characterizations are already known in the case when $X$ is finite, see \cite{TM}, \cite{T04}, \cite{T04s} from where we extract the following.
\begin{facts}
\label{qbsfinite}
 Suppose that $(X, r)$ is a finite quantum binomial set, $\mid X \mid = n,$
let $\Acal = \Acal(\textbf{k},X,r)$
be the associated quadratic algebra over a field $\textbf{k}$.

 Then any of the conditions (\ref{YBESquarefree}), ..., (\ref{YBElr3}) of Theorem \ref{qbs} is equivalent to any of following:
\begin{enumerate}
\item
\label{rsolutionofYBEb}
$S(X,r)$ is a semigroup of skew polynomial type (with respect to some appropriate enumeration of $X$).
\item
 \label{YBEASkew}
$\Acal$ is a binomial skew
 polynomial ring (with respect to some appropriate enumeration of $X$), that is an
 analogue of the Poincar\'{e}-Birkhoff-Witt
 theorem holds.
\item
 \label{YBEAfrobenius}
$(X,r)$ satisfies \textbf{lri} and $S(X,r)$ is a regular quantum monoid.
\item
 \label{YBESITYPE}
$S(X,r)$ is a semigroup of $I$-type, (see \cite{TM} for the definition).
\end{enumerate}
\end{facts}
Referring to condition (1) of the theorem briefly as YBE,  \cite[Theorem 1.3]{TM}
gives  (\ref{YBESITYPE}) $\Longleftrightarrow$ YBE,
\cite[Theorem 1.1]{TM} gives
(\ref{YBEASkew}) $\Longrightarrow$ YBE, while the
reverse implication (\ref{YBEASkew}) $\Longleftarrow$ YBE follows from \cite[Theorem 2.26]{T04}  (see
also \cite[Theorem B]{T04s}). Finally, (\ref{YBEAfrobenius})
 $\Longleftrightarrow$ YBE is shown in  \cite[Theorem B]{T04s}.

\section{Regular Garside monoids are  monoids of skew polynomial type}

In this section $(X,r)$ will denote a finite quantum binomial quadratic set with \textbf{lri},  $\mid X \mid = n.$
We shall assume that $S= S(x,r)$ is a regular Garside monoid, see Definition \ref{garsidedef}.

The following lemma will be used in the sequel. It is a
modification  of \cite[Lemma 2.19]{T04s}.

\begin{lemma}
\label{weakcyclicconditionL} Suppose $(X,r)$ is a finite quantum
binomial set with \textbf{lri}. Let $\Ocal=
\Ocal_{\Dcal_3}(\omega)$ be an arbitrary  orbit of the action of
$\Dcal_3$ on $X^3$, denote by $\Delta_i$ the diagonal of
$X^{\times i}, 2\leq i \leq 3$. Then the following conditions
hold.
\begin{enumerate}
\item $\Ocal \bigcap \Delta_3 \neq \emptyset$ if and only if
$\Ocal= \{xxx \},$ for some $x \in X$. \item $\Ocal \bigcap
((\Delta_2 \times X\bigcup X \times \Delta_2)\backslash \Delta_3))
\neq \emptyset$ if and only if $\mid \Ocal \mid = 3$. \item In
each of the cases   $\omega =yyx,$ or $\omega =yxx,$ where $x, y
\in X, x\neq y$, the orbit $\Ocal_{\Dcal_3}(\omega)$ contains
exactly $3$ elements. More precisely, the following implications
hold
\[\begin{array}{lcl}
(yx= x_1y_1) \in \Re & \Longrightarrow &( yx_1=x_2y_1)\in \Re,\; \text{and}\\
&&\Ocal_{\Dcal_3}(yyx)=\{ yyx, yx_1y_1,
x_2y_1y_1  \}, \\
(yx= x_1y_1) \in \Re &\Longrightarrow &(y_1x=x_1y_2) \in \Re,\; \text{and}\\
&&\Ocal_{\Dcal_3}(yxx)=\{ yxx, x_1y_1x, x_1x_1y_2 \},
\end{array}
\]
where $y, y_1, y_2$ respectively  $x, x_1,x_2$, are not
necessarily pairwise distinct.
\end{enumerate}
Furthermore, suppose $<$ is an ordering on $X$ such that every
relation in $\Re$ is of the type $yx=x^{\prime}y^{\prime},$ where
$y> x,$  $x^{\prime}< y^{\prime}$, and $y > x^{\prime}$. Then the
orbit $\Ocal_{\Dcal}(y_1y_2y_3)$ with $y_1\prec y_2 \prec y_3$
does not contain elements of the form $xxy,$ or $xyy$, $x\neq y
\in X$.
\end{lemma}

\begin{theorem}
\label{theoremA}
Let  $S=S(X,r)$ be a regular Garside monoid,
let
 $\Delta= x_1x_2...x_n$ be a regular presentation of its Garside element.
Then $S$ is a monoid of skew-polinomial type with respect to the canonical enumeration
$X= \{x_1, \cdots, x_n\}$ .
 \end{theorem}

 We shall verify the theorem in two steps.
 \textbf{Step I.} We prove that  every monomial $x_ix_j, 1 \leq i < j \leq n,$ is normal.

 \textbf{Step II.} We show that $\Re$ is a Gr\"obner basis, so condition (iii) of Definition \ref{skewpolynomialsemigroup} is also in force.

 All proofs are combinatorial and remind the proofs in \cite{T04s}. However, due to differences in hypothesis we
 write them explicitly.
 \begin{remark}
 \label{subwordsofdelta}
 Note that by hypothesis the monomial $\omega= x_1x_2\cdots x_n\in \langle X\rangle$ is in normal form and therefore each subword $a$ of $\omega$ is also in normal form, in particular each $x_ix_{i+1}, 1 \leq i \leq n-1,$ is in normal form.
 \end{remark}
 We shall use the terminology of \cite{T04s}.
\begin{definition}
\label{headsandtails} \cite{T04s} Let $S$ be a monoid, let $w \in
S.$ We say that $h \in X$ is \emph{a head of $w$}  if $w$ can be
presented (in $S$) as \[w = hw_1,\] where $w_1 \in \langle X
\rangle$ is a monomial of length $|w_1| = |w| -1$. Analogously, $t
\in X$ is \emph{a  tail of $w$} if \[w = w^{\prime}t\quad
\text{(in $S _0$)}\]  for some $w^{\prime}\in \langle X \rangle,$
with $|w^{\prime}| = |w| -1$.
\end{definition}
 \begin{notation}
 $H_{w}$ will denote the set of all heads of $w$ in $S$, respectively $T_w$
 will denote the set of all tails of $w$.
 \end{notation}
 \begin{lemma}
\label{Deltalemma} The Garside element $\Delta$
satisfies the conditions:
\begin{enumerate}
\item $\Delta$ is a monomial of length $n$.  There exist $n!$ distinct
words $\omega_i \in \langle X \rangle$,   $1 \leq i \leq n!$, for
which the equalities $\omega_i=\Delta$  hold in $S$. (These are the elements of the orbit
$\Ocal_{\Dcal_n}(x_1 \cdots x_n)$.  We call
them \emph{presentations} of $\Delta$. \item Every $x\in X$ is a left and right divisor of $\Delta$, so it occurs as a
``head''  (respectively, as a ``tail'')  of some presentation of
$\Delta.$ Thus there are equalities in $S:$
\[\begin{array}{cl}
\Delta&=x_1w_1^{\prime} =x_2w_2^{\prime}= \cdots x_nw_n^{\prime}\\
 &= \omega_1x_1=\omega_2x_2= \cdots \omega_nx_n.
 \end{array}
\]
 \item \label{W3} Every square-free monomial $a \in S$ of length
 $k, k \leq n$
 has exactly $k$ distinct ``heads'',
$h_1, \cdots, h_k$, and exactly $k$ distinct ``tails'', $t_1,
\cdots, t_k.$ The pair of sets $H_a, T_{a}$ uniquely determine $a$
in $S$. \item Furthermore for $1 \leq j < j+k \leq n$ the
monomials $\tau_{j j+k}= x_jx_{j+1}\cdots x_{j+k}$ satisfy
\[
H_{\tau_{j j+k}}= T_{\tau_{j j+k}}= \{x_j, x_{j+1}, \cdots, x_{j+k}\}
\]
\item \label{W4} $\Delta$ is the shortest monomial which
``encodes'' all the information about the relations $\Re.$
More precisely, for any relation $(xy=y^{\prime}x^{\prime}) \in
\Re,$ there exists an $a\in \langle X\rangle$, such that
$W_1=xya$ and $W_2=y^{\prime}x^{\prime}a$  are (different)
presentations of $W$.
\end{enumerate}
\end{lemma}


\begin{lemma}
\label{regular4}
For each integer $j,$ $1 \leq j \leq n-1,$ let
$\xi_{j,j+1}, \cdots,  \xi_{j,n}$, $\eta_{j,j+1}, \cdots,  \eta_{j,n}$
be the elements of $X$ uniquelly determined by the relations
\begin{equation}
\begin{aligned}
\label{xij} &(\xi_{j,j+1}\eta_{j,j+1}=x_jx_{j+1})\in \Re \cr
&(\xi_{j,j+2}\eta_{j,j+2}=\eta_{j,j+1}x_{j+2})\in \Re\cr &\cdots
\cdots \cdots\cr &(\xi_{j,n-1}\eta_{j,n-1}=\eta_{j,n-2}x_{n-1})
\in \Re\cr &(\xi_{j,n}\eta_{j,n}=\eta_{j, n-1}x_{n}) \in \Re. \cr
\end{aligned}
\end{equation}
Then for each $j$, $1 \leq j \leq n-1,$ the following conditions
are in force:
\begin{enumerate}
\item
\label{r1}
$\xi_{j,j+s}\neq \eta_{j,j+s-1},$ for all $s, 2\leq s \leq n-j$.
\item
\label{r2}
The equality
\[
\xi_{j,j+1}\xi_{j,j+2} \cdots \xi_{j,n} = x_{j+1} \cdots x_n\quad\text{holds in $S$}.
\]
\item
\label{r3}
$x_{j+1}x_{j+2}\cdots x_n \eta_{j,n}= x_jx_{j+1}\cdots x_n$ holds in $S$.
\item
\label{r4}
The elements $\eta_{j,n}, \eta_{j+1,n}, \cdots, \eta_{n-1,n}$ are
pairwise distinct.
\end{enumerate}
\end{lemma}
\begin{proof}
Condition (\ref{r1}) follows from the Ore conditions on generators, see (2) and (3) of Proposition
\ref{nondeg_and_invol_goodprop}.
To prove the remaining conditions
we use decreasing induction on $j$, $1 \leq j \leq n-1.$

\textbf{Step 1.} $j=n-1.$ Clearly, $x_{n-1}x_n$ is in normal form,
so the relation in $\Re $ in which it occurs has the shape
$x_{n-1}x_n=\xi_{n-1,n}\eta_{n-1,n},$ with $\xi_{n-1,n}> x_{n-1}.$
It follows then that $\xi_{n-1,n}=x_n$ and
$x_{n-1}x_n=x_n\eta_{n-1,n}.$
This gives (\ref{r2}) and (\ref{r3}). There is nothing to prove  in  (\ref{r4}).

\textbf{Step 2.} We first prove (\ref{r4}) for all $j, 1 \leq j<
n-1.$ Assume that for all $k, \; n-1\geq k > j,$ the elements
$x_k, x_{k+1},\cdots, x_n$, $\xi_{k,k+1}, \cdots, \xi_{k,n}$,
$\eta_{k,k+1},\cdots, \eta_{k,n}$ satisfy
\begin{equation}
\begin{aligned}
\label{xik} &\xi_{k,k+1}\eta_{k,k+1}=x_kx_{k+1}\in \Re \cr
&\xi_{k,k+2}\eta_{k,k+2}=\eta_{k,k+1}x_{k+2}\in \Re  \cr &\cdots
\cdots \cdots\cr &\xi_{k,n-1}\eta_{k,n-1}=\eta_{k,n-2}x_{n-1} \in
\Re \cr &\xi_{k,n}\eta_{k,n}=\eta_{k, n-1}x_{n} \in \Re , \cr
\end{aligned}
\end{equation}
all $\eta_{j+1,n}, \eta_{j+2,n}, \cdots,\eta_{n-1,n}$
are pairwise distinct, and the modified
conditions
(\ref{r4}), in which ``$j$'' is replaced by ``$k$'' hold. Let
$\xi_{j,j+1}, \cdots \xi_{j,n}$, $\eta_{j,j+1}\cdots \eta_{j,n}$
satisfy (\ref{xik}). We shall prove  that $\eta_{j,n}\neq
\eta_{k,n},$ for all $k, j<k\neq n-1.$ Assume the contrary,
\[
\eta_{j,n}=\eta_{k,n}
\]
for some $k>j.$
Consider the relations
\begin{equation}
\label{e2} \xi_{j,n}\eta_{j,n}=\eta_{j,n-1}x_n, \quad
\xi_{k,n}\eta_{k,n}=\eta_{k,n-1}x_n.
\end{equation}
 By Proposition
\ref{nondeg_and_invol_goodprop} the Ore condition holds, so
(\ref{e2}) imply
\[
\label{e3} \eta_{j,n-1}=\eta_{k,n-1}.
\]
Using the same argument in $n-k$ steps we
obtain the equalities
\[
\eta_{j,n}=\eta_{k,n}, \eta_{j,n-1}=\eta_{k,n-1}, \cdots,  \eta_{j,k+1}=\eta_{k,k+1}.
\]
Now the relations
\[
\xi_{j,k+1}\eta_{j,k+1}=\eta_{j,k}x_{k+1},  \;
\xi_{k,k+1}\eta_{k,k+1}=x_kx_{k+1},
\]
and the Ore condition  imply $\eta_{j,k}= x_k.$
Thus, by (\ref{xik}) and (\ref{xij}) we obtain  a relation
\[\xi_{j,k}x_k=\xi_{j,k-1}x_k \in \Re.\] This is impossible, by
Proposition \ref{nondeg_and_invol_goodprop}. We have shown that the
assumption $\eta_{j,n}=\eta_{k,n}$, for some $k>j$, leads to a
contradiction. This proves (\ref{r4})
for all $j, 1 \leq j\leq n-1$.

We set
\begin{equation}
\label{eta}
\eta_1=\eta_{1,n}, \eta_2=\eta_{2,n}, \cdots, \eta_{n-1}=\eta_{n-1,n}.
\end{equation}

Next  we prove (\ref{r2}) and (\ref{r3}).

By the inductive assumption, for $k > j,$ we have
\[
\label{e4}
\xi_{k,k+1}\xi_{k,k+2}\cdots\xi_{k,n} = x_{k+1} \cdots x_n,
\]
and
\[
x_{k+1} \cdots x_n.\eta_{k+1}=x_k\cdots x_n.
\]
Applying the relations (\ref{xik}) one easily sees,  that
\[
\xi_{j,j+1}\xi_{j,j+2}\cdots\xi_{j,n}.\eta_{j,n}=x_jx_{j+1}\cdots x_n.
\]
Denote
\[
\omega_j=\xi_{j,j+1}\xi_{j,j+2}\cdots \xi_{j,n}.
\]
We have to show that the normal form,  $Nor(\omega_j)$, of
$\omega_j$ satisfies the equality of words in $\langle X \rangle$
\[Nor(\omega_j) =
x_{j+1}x_{j+2} \cdots x_n.\] As a subword of length $n-j$ of the
presentation $\Delta =x_1x_2\cdots x_{j-1}w_j\eta_{j,n}$, the
monomial  $\omega_j$ has exactly $n-j$ heads
\begin{equation}
\label{e06} h_1 < h_2 < \cdots < h_{n-j}.
\end{equation}
But $Nor(\omega_j)=\omega_j$, is an equality  in $\Scal$, so the
monomial $Nor(\omega_j)$ has the same heads as $\omega_j.$
Furthermore, there is an equality of words in $\langle X \rangle$,
  $Nor(\omega_j)=h_1\omega^{\prime}$, where $\omega^{\prime}$ is a
monomial of length $n- j -1.$ First we see that $h_1  \geq x_{j}$
This follows immediately from the properties of the normal
monomials and the relations
\begin{equation}
\label{e7}
\Nor(\omega_j)\eta_j=\omega_j\eta_j= x_jx_{j+1}\cdots x_n \in N .
\end{equation}
Next we claim that $h_1 > x_j.$
Assume the contrary, $h_1=x_j.$
Then by (\ref{e7}) one has
\[
x_j\omega^{\prime}\eta_j=\omega_j\eta_j=x_jx_{j+1}\cdots x_n.
\]
As a Garside monoid, $S$ is cancellative, hence
\[
\omega^{\prime}\eta_j=x_{j+1}\cdots x_n \in N.
\]
Thus $\eta_j$ is a tail of the monomial $x_{j+1}\cdots x_n.$ By
the inductive assumption, conditions (\ref{r2}) and (\ref{r3}) are
satisfied, which together with (\ref{eta}) give additional $n-j$
distinct tails of the monomial $x_{j+1}\cdots x_n$, namely
$\eta_{j+1}, \eta_{j+2}, \cdots \eta_{n-1}, x_n$. It follows then
that the monomial $x_{j+1}\cdots x_n$ of length $n-j$ has $n-j+1$
distinct tails, which is impossible. This implies  $h_1 > x_j.$
Now since $\omega_j$ has precisely $n- j = |\omega_j|$ distinct
heads, which in addition satisfy (\ref{e06}) we obtain equality of
sets
\[
H_{\omega_j}=\{h_1, h_2, \cdots , h_{n-j}\}=\{x_{j+1}, x_{j+2}, \cdots,  x_n\}.
\]
Recall that by the inductive assumption we have
\begin{equation}
\label{indasseq}
H_{x_{j+1}x_{j+2}\cdots x_n} = \{x_{j+1}, x_{j+2}, \cdots,  x_n\}.
\end{equation}
The  equality $H_{\omega_j}= H_{x_{j+1}x_{j+2}\cdots x_n}$ clearly
implies equality of monomials in $S:$
\[\omega_j=
x_{j+1}x_{j+2}\cdots x_n.  \] We have shown
(\ref{r3}). Now the equalities
\[
x_{j+1}\cdots x_n\eta_j=x_{j}x_{j+1}\cdots x_n
\]
and (\ref{indasseq}) yield :
\[H_{x_{j}x_{j+1}\cdots x_n} = \{x_{j}, x_{j+1}, x_{j+2}, \cdots,  x_n\},\]
which  proves (\ref{r2}). The lemma has been proved.
\end{proof}

\begin{proposition}
\label{GarsideimpliesSkewtypeProp} In notation as in Lemma \ref{regular4} the following conditions hold in $S$.
\begin{enumerate}
\item \label{p1} For any integer $j,$ $1 \leq j \leq n-1,$ there
exists unique $\eta_j \in X,$ such that
\[
x_{j+1} \cdots x_n\eta_j = x_jx_{j+1}\cdots x_n.
\]
\item \label{p3} The elements $\eta_1, \eta_2, \cdots, \eta_{n-1}$
are pairwise distinct. \item \label{p4} For each $j$, $1 \leq j
\leq n-1,$ the set of heads $H_{W_j}$ of the monomial $W_j =
x_jx_{j+1}\cdots x_n$ is
\[
H_{W_j}= \{ x_j, x_{j+1},  \cdots , x_n \}.
\]
\item
\label{p5}
For any pair of integers $i,j$, $1 \leq i < j \leq n,$
the monomial $x_ix_j$ is normal. Furthermore, the unique
relation in which $x_ix_j$ occurs
has the form $x_{j^{\prime}}x_{i^{\prime}}= x_ix_j$, with
$j^{\prime} > i^{\prime},$ and $j^{\prime} > i.$
\end{enumerate}
\end{proposition}

\begin{proof}
Conditions  (\ref{p1}), (\ref{p3}), (\ref{p4}) of the proposition follow from
Lemma \ref{regular4}, here $\eta_j$ are defined in (\ref{eta}) . We shall prove first that for any pair
$i,j,$ $1\leq i < j \leq n,$ the monomial $x_ix_j$ is normal.
Assume the contrary. Then there is a relation
\begin{equation}
\label{1p}
(x_ix_j=x_{j^{\prime}}x_{i^{\prime}}) \in \Re,\quad\text{where}\quad j^{\prime}< i.
\end{equation}
Clearly,
\begin{equation}
\label{equ1} (x_jx_{j+1}\cdots x_n)\eta_{j-1} \eta_{j-2} \cdots
\eta_{i+1}= x_{i+1}x_{i+2} \cdots x_n,
\end{equation}
holds in $S$. Consider the monomial $u = x_ix_jx_{j+1}\cdots
x_n\eta_{j-1} \eta_{j-2} \cdots \eta_{i+1}\in X^{n-i+1}.$
 It satisfies the following equalities in $S$
 \[
 \begin{array}{ll}
 u  &= x_ix_jx_{j+1}x_{j+2}\cdots x_n\eta_{j-1} \eta_{j-2} \cdots \eta_{i+1} \\
    &= x_{i+1}x_{i+2} \cdots x_n\quad \text{by (\ref{equ1})}\\
    & = x_{j^{\prime}}x_{i^{\prime}} x_{j+2}\cdots x_n\eta_{j-1} \eta_{j-2} \cdots \eta_{i+1}                            \quad \text{by (\ref{1p})}
 \end{array}
 \]
 Note that the monomial $x_{i+1}x_{i+2} \cdots x_n$ is in normal form, therefore $\Nor u = x_{i+1}x_{i+2} \cdots x_n$
 The inequality $j^{\prime}< i$ implies the obvious inequality in $\langle X\rangle:$
 \[x_{j^{\prime}}x_{i^{\prime}} x_{j+2}\cdots x_n\eta_{j-1} \eta_{j-2} \cdots \eta_{i+1} < x_{i+1}x_{i+2} \cdots x_n = \Nor u,\]
which is impossible, since $\Nor u$ is the minimal element in the
orbit $\Dcal_{n-i+1}(u).$ We have proved that all monomials
$x_ix_j$ with $1\leq i < j \leq n$ are in normal form. Note that
the number of relations is exactly $\frac{n(n-1)}{2}$ and each
relation contains exactly one normal monomial, it follows then
that  each monomial $x_jx_i,$ with  $n\geq j > i\geq 1$, is not in
normal form. Hence,
 each relation in $\Re$ has the shape
$y_jy_i=y_{i^{\prime}}y_{j^{\prime}},$ where $1 \leq i < j \leq
n,$ $1 \leq i^{\prime} < j^{\prime} \leq n,$ and $j > i^{\prime}$,
which proves (\ref{p4}) and (\ref{p5}).
\end{proof}
Step I is complete now. This straightforwardly implies that
 the (unique) relation in $\Re$ in which  $x_ix_j$ occurs has the shape
 \[x_{j^{\prime}}x_{i^{\prime}}= x_ix_j, \quad\text{with}\quad j^{\prime} > i^{\prime},  j^{\prime} > i, \]
 so the relations $\Re$ have the ``correct shape"  \emph{of skew-polynomial type}, thus conditions (\ref{skewrelations}) (i) and (ii)
 of Definition \ref{skewpolynomialsemigroup} hold.

   We next proceed with \textbf{Step II}.

\begin{lemma}
\label{GarsideimpliesSkewtypeLemma}
\begin{enumerate}
\item[(a)] \label{GarsideimpliesSkewtypeLemmaa} The set of
relations $\Re $  is a Gr\"{o}bner basis with respect to the
ordering $< $ on $\langle X \rangle.$ \item[(b)]
\label{GarsideimpliesSkewtypeLemmab} $S(X,r)= \langle x_1 x_2
\cdots x_n \mid \Re \rangle$ is a semigroup of skew polynomial
type.
\end{enumerate}
\end{lemma}

\begin{proof}
It will be enough to prove that the ambiguities $y_ky_jy_i,$ with  $k > j >i$,
do not give rise to new relations
in $S.$ Or equivalently, each ordered monomial $\omega$ of length $3$,
that is $\omega = xyz, x\leq y\leq z,$  is in normal form.

Let $\omega$ be an ordered monomial of length $3$, and assume
$Nor_{\omega} \neq \omega$ (as words in $\langle X \rangle$). This
means there exists an $\omega^{\prime}\in \Dcal_3(\omega)$, such
that $\omega^{\prime} < \omega$ (in the degree-lexicographic
ordering $<$ in $\langle X \rangle$).  Four cases are possible:
\begin{equation}\begin{array}{ll}
\label{N1}
(i)\quad &\omega = x_ix_jx_k, 1 \leq i < j < k \leq n\\
(ii)\quad &\omega = x_ix_ix_j, 1 \leq i < j  \leq n\\
(iii)\quad &\omega = x_ix_jx_j, 1 \leq i < j  \leq n\\
(iv) \quad &\omega = x_ix_ix_i, 1 \leq i  \leq n.
 \end{array}
\end{equation}
Suppose (\ref{N1}) \textbf{(i)} holds. Then there is an equality
of elements in $S:$
\[\omega = x_ix_jx_k =
x_{i^{\prime}}x_{j^{\prime}}x_{k^{\prime}} =\omega^{\prime} ,\; \text{where} \; x_{i^{\prime}}\leq
x_{j^{\prime}}\leq x_{k^{\prime}},\] and  as elements of $\langle X
\rangle$, the two monomials satisfy
\begin{equation}
\label{N5}  \omega^{\prime}= x_{i^{\prime}}x_{j^{\prime}}x_{k^{\prime}}< x_ix_jx_k = \omega.
\end{equation}
By (\ref{N5}) one has
\begin{equation}
\label{N6}  x_{i^{\prime}}\leq x_i .
\end{equation}
We claim that there is an inequality $x_i^{\prime} < x_i.$ Indeed,
it follows from Lemma \ref{weakcyclicconditionL} that the orbit
$\Ocal_{\Dcal_3}(y_iy_jy_k)$ does not contain elements of the
shape $xxy,$ or $xyy,$ therefore an assumption,
$x_i=x_{i^{\prime}}$ would imply
$x_jx_k=x_{j^{\prime}}x_{k^{\prime}}$ with $x_j < x_k$ and
$x_{j^{\prime}}< x_{k^{\prime}}$, which contradicts Proposition
\ref{GarsideimpliesSkewtypeProp}. We have obtained that
$x_{i^{\prime}}\leq x_i$. One can easily see that there exists an
$\omega \in \langle X \rangle,$ such that
\[
(x_ix_jx_k) * \omega = x_i x_{i+1}\cdots x_n .
\]
The monomial $x_i x_{i+1}\cdots x_n$ is in normal form, hence
\[\Nor ((x_ix_jx_k) * \omega) = x_i x_{i+1}\cdots x_n\].  But
$x_{i^{\prime}}x_{j^{\prime}}x_{k^{\prime}} = x_ix_jx_k$ in $S$,
hence
\[
Nor(x_{i^{\prime}}x_{j^{\prime}}x_{k^{\prime}} * \omega) = Nor (x_ix_jx_k * \omega ) =x_i x_{i+1}\cdots x_n.
\]
This together with the following inequalities in $\langle X \rangle:$
\[
Nor(x_{i^{\prime}}x_{j^{\prime}}y_{k^{\prime}} * \omega) \leq x_{i^{\prime}}x_{j^{\prime}}x_{k^{\prime}} * \omega < x_i x_{i+1}\cdots x_n
\]
give a contradiction. It follows then that monomial $x_ix_jx_k,$
$i < j < k,$ is in normal form.

Suppose now \textbf{(ii)} $\omega = x_ix_ix_j, 1 \leq i < j  \leq
n$. It follows from Lemma \ref{weakcyclicconditionL}  that the
orbit $\Ocal = \Ocal_{\Dcal_3}(x_ix_ix_k)$ is the set
\[
\Ocal= \{ \omega= x_ix_ix_j,\;\; \omega_2=
x_ix_{j^{\prime}}x_{i^{\prime}},\;\;
\omega_3=x_{j^{\prime\prime}}x_{i^{\prime}}x_{i^{\prime}} \}
\]
where
\[
x_{j^{\prime}}x_{i^{\prime}}=x_ix_j \in \Re, \; \text{with} \; x_i<
x_{j^{\prime}}> x_{i^{\prime}}
\]
and
\[
x_{j^{\prime\prime}}x_{i^{\prime}}=x_ix_{j^{\prime}}\in \Re, \;
\text{with} \; x_{j^{\prime\prime}}> x_{i^{\prime}}.
\]
It is clear then that there are strict inequalities in $\langle X\rangle:$
\[
\omega= x_ix_ix_j < \omega_2= x_ix_{j^{\prime}}x_{i^{\prime}} <
\omega_3=x_{j^{\prime\prime}}x_{i^{\prime}}x_{i^{\prime}}
\]
which gives $Nor(\omega) = \omega.$

The case \textbf{(iii)}  is analogous to (ii). Case \textbf{(iv)}
is straightforward, since all relations are square free. We have
proved that all ordered monomials of length $3$ are in normal
form, and therefore the ambiguities $x_kx_jx_i$ do not give rise
to new relations in $S$ (i.e. are solvable). It follows then that
the set of relations $\Re$ is a Gr\"{o}bner basis.

We have verified all conditions in  Definition \ref{skewpolynomialsemigroup}, therefore  $S(X,r)$ is a monomial of skew polynomial type. The lemma has been proved.
\end{proof}
\begin{prooftheorem}
Now Theorem \ref{theoremA}
follows straightforwardly from Proposition  \ref{GarsideimpliesSkewtypeProp} and Lemma \ref{GarsideimpliesSkewtypeLemma}.
\end{prooftheorem}

\section{Proof if the main theorem}
\label{section_regularity}
The Koszul dual $A^{!}$ of a quadratic algebra $A$ was introduced by Yu. Manin, see   \cite{Maninpreprint}, 5.1. The properties of the two algebras $A$ and $A^{!}$ are closely related.

Suppose $(X,r)$ is a finite quantum binomial quadratic set, for convenience we enumerate $X = \{x_1, \cdots x_n\}$.
(In most of the cases the enumeration will be induced by a regular monomial of length $n$).
Suppose $\Re_0=\Re_0(r)$, see (\ref{Re0}) is the set of defining relations of the quantum binomial algebra
$\Acal= \Acal(\textbf{k},X,r)$, so  $\Acal= \textbf{k}\langle X \rangle/(\Re_0)$. Clearly, $\Acal$ is a quadratic algebra, and one can extract from \cite{Maninpreprint}
an explicit presentation of its Koszul dual, $\Acal^{!}$.
\begin{definition}
\label{koszuldualdef} \emph{The Koszul dual $\Acal^{!}$
of $\Acal$},   is the quadratic algebra,
\[
\Acal^{!}:=
 \textbf{k}\langle \xi _1, \cdots, \xi _n \rangle/(\Re_0^{\bot}),
\]
where the set $\Re_0^{\bot}$ contains precisely  $\binom{n}{2}+n$
relations  of the following two types:

a) \emph{binomials:}
\[
\xi_j\xi_i+\xi_{i^{\prime}}\xi_{j^{\prime}}\in
\Re^{\bot}, \ \text{whenever} \
x_jx_i-x_{i^{\prime}}x_{j^{\prime}}\in \Re_0, \ 1\leq i\neq j
\leq n; \;
\]
and
b) \emph{monomials:}
\[
(\xi_i)^2 \in \Re^{\bot} , 1 \leq i \leq n.
\]
\end{definition}
\begin{remark} \cite{Maninpreprint}
Consider the vector spaces   $V=\Span_{\textbf{k}}(x_1, x_2,
\cdots, x_n),$  $V^{*}= \Span_{\textbf{k}}(\xi_1 ,\xi_2 \cdots,
\xi_n),$ and define a bilinear pairing
\[\langle \; \mid \; \rangle : V^{*}\otimes
V\longrightarrow \textbf{k}\quad\text{by}\quad \langle \xi _i\mid
x_j\rangle : = \delta_{ij}.
\]
 Then the relations $\Re_0^{\bot}$ generate a subspace
in $V^{*}\otimes V^{*}$ which is orthogonal to the subspace of $
V\otimes V$ generated by $\Re_0.$
\end{remark}
\begin{definition}
\label{S!def} We introduce a monoid $S^{!}= S^{!}(X,r)$  with
\emph{zero element} denoted by $0$. $S^{!}$ is generated by $X$
and has a presentation as
\[
 S^{!} := \langle X \mid  \Re^{!}\rangle
 \]
where the set $ \Re^{!}$ contains  $\binom{n}{2}+n$
defining relations
\[
 \Re^{!} = \Re^{!}(r) := \Re (r) \bigcup \{xx=0\mid x \in X \}
 \]
\end{definition}
Clearly, there is a canonical epimorphism $\varphi: S(X,r) \longrightarrow S^{!}(X,r),$
\[
\varphi (u) = \left\{ \begin{array}{l} u,\; \text{if}\; u \in S \; \text{is a square-free monomial}\\
                                      0 \quad \text{else}.
\end{array}\right\}
\]
\begin{proposition}
\cite{T04s}
Suppose  we can fix an enumeration $X = \{x_1, \cdots, x_n \}$,
such that each relation in $\Re_0$ has the shape
\[x_jx_i-x_{i^{\prime}}x_{j^{\prime}},\quad
1 \leq i < j \leq n, \;1 \leq i^{\prime} < j^{\prime} < n, \;\; j
> i^{\prime}.\] As usual, we shall consider Gr\"{o}bner bases with
respect to degree -lexicographic  ordering  on $\langle X\rangle$.
The  following are equivalent.
\begin{enumerate}
\item
The set $\Re_0$ is a Gr\"{o}bner basis of the ideal $(\Re_0),$ or equivalently, the algebra $\Acal = \Acal(\textbf{k},X,r)$ is a binomial skew polynomial ring.
\item
$S(X,r)$ is a monoid of skew polynomial type.
\item
The set of ordered monomials
 \[
\Ncal= \{x_1^{\alpha_1}x_2^{\alpha _2}\cdots x_n^{\alpha_n} \mid 0
\leq \alpha _i \, \ \text{for} \ 1 \leq i \leq n \}
\]
is a $\textbf{k}$-basis of $\Acal.$
\item
The set $\Re^{\bot}$ is a Gr\"{o}bner basis of the ideal $(\Re^{\bot}).$
\item
The set $(\Re^{!})_0$ is a Gr\"{o}bner basis,  of the ideal $((\Re^{!})_0).$
\item
The set of monomials
\[
\Ncal^{!}= \{x_1^{\varepsilon_1}x_2^{\varepsilon _2}\cdots
x_n^{\varepsilon_n} \mid 0 \leq\varepsilon _i \leq 1, \ \text{for}
\ 1 \leq i \leq n \}
\]
is a $\textbf{k}$-basis of the Koszul dual algebra $\Acal^{!}$.
\item
$\Ncal^{!}$ is a $\textbf{k}$-basis of the semigroup algebra $\textbf{k} S^{!}$,
so  $S^{!}$ can be identified as a set with $\Ncal^{!}$.
\end{enumerate}
\end{proposition}
\begin{definition}
\label{frobeniusalgebra}\cite{Maninpreprint}, \cite{Maninbook} A
graded algebra $A=\bigoplus_{i\ge0}A_i$ is called
\emph{a Frobenius algebra of dimension $n$},
(or \emph{a Frobenius quantum space of dimension $n$})
if
\begin{enumerate}
\item[(a)]
$dim(A_n)=1$,  $A_i =0,$ for $i>n.$
\item[(b)]
For all $n\geq j \geq 0$) the multiplicative map
$m: A_j\otimes A_{n -j} \rightarrow A_n$
is a perfect duality (nondegenerate pairing).

A Frobenius algebra $A$ is called \emph{a quantum Grassmann algebra} if in addition
\item[(c)]
$dim_ \textbf{k}A_i= \binom{n}{i} , \;\text{for}\; 1 \leq i \leq
n.$ \end{enumerate}
\end{definition}
\begin{definition}
\label{socleandprincipalmonomialdef}\cite{T04s}
In notation as above  suppose the  Koszul dual
  $A^{!} = A^{!}(\Re_0)$ is Frobenius.
  The one dimensional component $A_n^{!}$ is called \emph{the socle
of} $A^{!}.$  Let $\Wcal=\Wcal(\xi_1, \xi_2, \cdots, \xi_n)$ be the monomial which spans  the socle.
In this case the monoid $S^{!}$ has top degree $n$ and contains exactly one (nonzero) element $W=W(x_1, \cdots, x_n)$ of length $n$.
Clearly, its pull-back $\varphi^{-1}(W)$ in $S$, also contains exactly one element,  $W$.
 Note that in this case $W$ is the unique square-free element  in $S$ with  length $|W|=n$. We call it \emph{the principal monomial of $S$}
 (respectively \emph{the principal monomial of $\Acal$}), see \cite{T04s}.
 Suppose  $W$ can be presented as a product of exactly $n$ generators,
  $W = y_1 \cdots, y_n, y_i \in X,$ $y_i\neq y_j$ whenever $i \neq j$. So $X= \{y_1, \cdots y_n\}.$ Fix the
  degree-lexicographic ordering $<$ on $\langle X\rangle$, where $y_1 < \cdots < y_n$.
 Denote $\omega_0 = y_1 \cdots y_n$ (considered as a word in $\langle X\rangle$), and look at the orbit
  $\Ocal= \Ocal_{\Dcal_n}(\omega_0)$,  in $X^k$. ( $\omega\in \Ocal$ \emph{iff} $\omega = \omega_0= W$ as elements of $S$).
  By Definition
  \ref{regquantummondefinitions},
 $W = y_1 \cdots y_n$ is \emph{a regular presentation} of $W$ if  $\omega_0 $  is the minimal element of $\Ocal$ with respect to the degree-lexicographic ordering $<$.
  The principal monomial $W$ is \emph{regular} if it has some regular presentation.
  In this case  $\Acal^{!}$ is said to have \emph{a regular socle}.
 \end{definition}
 The following can be extracted from
 \cite[Theorem B]{T04s} (bearing in mind that lri implies the  cyclic conditions, see Remark \ref{lriremark2}).
 \begin{facts}\cite{T04s}
 \label{factstheoremB}
 Let $(X,r)$ be a quantum binomial set with \textbf{lri}.
$S = S(X,r),  \; \Acal= \Acal(\textbf{k}, X,r) = \textbf{k}\langle X \rangle/(\Re_0(r))$, as usual.
 Then the following three conditions
are equivalent.
\begin{enumerate}
\item \label{theoremB1} $S(X,r)$ is a regular quantum monoid, that
is the Koszul dual $A^{!}$ is Frobenius of dimension $n$, and  has
a regular socle.
\item
\label{theoremB2}
$S$ is a binomial skew
polynomial  monoid, with respect to some appropriate enumeration
of $X.$
\item
\label{theoremB2Alg} $\Acal(\textbf{k}, X,r)$ is a binomial skew
polynomial ring, with respect to some appropriate enumeration of
$X.$ \item \label{theoremB3} $(X,r)$ is a solution of YBE.
\end{enumerate}
\end{facts}
The following result shows that a milder assumption: $A^{!}(X,r)$
is Frobenius, and $S$ is cancellative, ($W$ is not necessarily
regular),  imply that $(S,W)$ is a Garside monoid.

\begin{theorem}
\label{frobeniusimpliespregarsidethm}
Let $(X,r)$ be quantum binomial set with \textbf{lri}. Suppose the Koszul dual $A^{!}(X,r)$ is Frobenius.
Let $W$ be the principal monomial in the monoid $S= S(X,r).$ Then
\begin{enumerate}
\item The set $\;\sum_{W}\;$ of all left divisors of $W$ in $S$
consists exactly of all square-free elements $u \in S$ of length
$|u|  \leq n$ and coincides with the set of all right divisors of
$W$. \item Furthermore, if $S$ is cancellative, then $W$ is a
Garside element, and $(S, W)$ is a comprehensive Garside monoid.
\end{enumerate}
\end{theorem}
\begin{proof}
Note that $W$ is square-free element of $S$  with length $|W| =
n$, see Definition \ref{socleandprincipalmonomialdef}. Then,
clearly, each left (and right) divisor $a$ of $W$ has length $\leq
n$ and is square-free. Conversely, let $a\in S$ be a square-free
element. It follows from the Frobenius property that $a$ has
length $\leq n$, and if $a$ has length $n$ then $a=W$. Assume
$m=|a| <n$. One has $0 \neq a^{\prime}=\varphi(a) \in
\Acal^{!}_m$, so by the Frobenius property, see Definition
\ref{frobeniusalgebra} (b), there exists a monomial $b \in
\Acal^{!}_{n-m}$ such that $0 \neq ab\in \Span W = \Acal^{!}_n.$
Clearly then $ab = W$ holds in $S$, that is $a$ is a left divisor
of $W$. Analogously one shows that $a$ is a right divisor of $W.$
We have verified that each square-free element $a \in S$ (of
length $\leq n$) is a left and a right divisor of $W$.

Assume furthermore that $S$ is cancellative. Then condition (1. b) of Definition \ref{garsidedef} is satisfied, so $W$
is a Garside element and (S, W) is a pre-Garside monoid. We claim that $(S, W)$ is a Garside monoid, that is condition (2) in
Definition \ref{garsidedef} holds.  $S$ is atomic with set of atoms $X$, see Remark \ref{atomicremark}.
 By hypothesis $(X,r)$ is a quantum binomial set with \textbf{lri} so by Proposition \ref{nondeg_and_invol_goodprop} the left and right Ore conditions on generators are satisfied. Hence for arbitrary pair $s,t \in X, s \neq t,$ there are uniquely determined $x,y \in X$ such that $sx=ty$ holds in $S$, the   element $\Delta_{s,t}= sx = ty$ is square-free and therefore is a left divisor of $W$. Clearly $\Delta_{s,t}$ is the unique minimal element of the set
\[\{a \in \sum\mid s \preceq a,\;\; \text{and}\;\; t \preceq a\}\]
with respect to the partial ordering $\preceq$ on $S$ induced by the left divisibility, ($u \preceq v$ iff $u$ is a left divisor of $v$). It is straightforward that $(S, W)$ is comprehensive, see Definition \ref{garsidedef} (3).
\end{proof}
\begin{corollary}
\label{theoremC}
Suppose $S= S(X,r)$ is a monoid of skew-polynomial type, with respect to the enumeration
$X = \{x_1, \cdots, x_n\}$. Then $S$ it is a regular Garside monoid with regular Garside element
$\Delta = x_1 x_2 \cdots x_n.$
 \end{corollary}
\begin{proof}
 By \cite[Theorem 5.13]{T96}
  $S$ is cancellative. Furthermore, the Koszul dual $\Acal^{!}$ is Frobenius, see \cite[Theorem 3.1]{T04s} ,  and $W = x_1 \cdots x_n$ is a regular presentation of the principal monomial of $S$.
It follows then by  Theorem \ref{frobeniusimpliespregarsidethm} that $(S, W)$ is a Garside monoid, and clearly its  Garside element $W$ has regular presentation.
\end{proof}
\begin{proofmaintheorem} Let $(X,r)$ be a quantum binomial set with \textbf{lri}. For the conditions (1), (2), (3) of the theorem we have the following implications:
\[\begin{array}{lll}
(1) &\Longrightarrow& (2) \quad \text{by Theorem \ref{theoremA}}\\
(2) &\Longleftrightarrow& (3)\quad \text{by \cite[Theorem B]{T04s},  see also Facts \ref{factstheoremB}}\\
(3) &\Longrightarrow & (1) \quad \text{by Corollary \ref{theoremC}}
\end{array}
\]
By Theorem \ref{frobeniusimpliespregarsidethm}, $S$ is a comprehensive Garside monoid with  Garside element $\Delta$. By \cite{KasselT},  the monoid
$S(X,r)$ is embedded in its group of fraction, which is isomorphic to the group $G= G(X,r)$. Thus  $(G,\Delta)$ is a Garside group.
\end{proofmaintheorem}

{\bf Acknowledgments}.  This paper was written during my visit to
Abdus Salam International Centre for Theoretical Physics (ICTP),
Trieste, Summer 2009. It is my pleasant duty to thank  Professor
Le Dung Trang and  the Mathematics group of ICTP for inviting me
and for our valuable and stimulating discussions. I thank ICTP for
the continuous support and for the inspiring and creative
atmosphere. Some ideas in this paper were inspired while giving a
course on Braid groups in the University of Granada, 2009. I thank
Pascual Jara for the invitation, for our interesting discussions
and for his amazing hospitality.

\bibliography{mybibs}
\bibliographystyle{amsabbrv}
\ifx\undefined\bysame
\newcommand{\bysame}{\leavevmode\hbox to3em{\hrulefill}\,}
\fi

\end{document}